# Unsupervised Surrogate-Assisted Synthesis of Free-Form Planar Antenna Topologies for IoT Applications

Khadijeh Askaripour[1], Adrian Bekasiewicz[1], and Slawomir Koziel[1,2]

[1]Faculty of Electronics, Telecommunications and Informatics, Gdansk University of Technology, Gdansk, Poland

[2] Department of Engineering, Reykjavik University, Reykjavik, Iceland

**Abstract**

Design of antenna structures for Internet of Things (IoT) applications is a challenging problem. Contemporary radiators are often subject to a number of electric and/or radiation-related requirements, but also constraints imposed by specifics of IoT systems and/or intended operational environments. Conventional approaches to antenna design typically involve manual development of topology intertwined with its tuning. Although proved useful, the approach is prone to errors and engineering bias. Alternatively, geometries can be generated and optimized without supervision of the designer. The process can be controlled by suitable algorithms to determine and then adjust the antenna geometry according to the specifications. Unfortunately, automatic design of IoT radiators is associated with challenges such as determination of desirable geometries or high optimization cost. In this work, a variable-fidelity framework for performance-oriented development of free-form antennas represented using the generic simulation models is proposed. The method employs a surrogate-assisted classifier capable of identifying a suitable radiator topology from a set of automatically generated (and stored for potential re-use) candidate designs. The obtained geometry is then subject to a bi-stage tuning performed using a gradient-based optimization engine. The presented framework is demonstrated based on six numerical experiments concerning unsupervised development of bandwidth-enhanced patch antennas dedicated to work within 5 GHz to 6 GHz and 6 GHz to 7 GHz bands, respectively. Extensive benchmarks of the method, as well as the generated topologies are also performed.

## 1. Introduction

High-performance of antenna structures is one of the key requirements behind the development of agile Internet of Things (IoT) systems dedicated to operate in dynamic environments. The electrical and field-related characteristics of radiators stem from the selected topology, architecture of the feeding network (especially for antenna arrays), as well as appropriate tuning of the design parameters [1]-[8]. The main problem with determination of antenna geometry revolves around the inverse flow of information [10]. In other words, the designer knows the specification requirements that need to be addressed, yet lacks the details on topologies and/or their modifications (other than the ones introduced based on experience) that can produce satisfactory responses [9]-[11].

Antennas for contemporary IoT systems are often subject to stringent requirements that might include broadband, multi-band, or otherwise unconventional behavior [1]-[4]. At the same time, the inferior performance of conventional radiators (e.g., synthesized using well-known empirical formulas) makes them unsuitable for modern applications [12], [13]. A practical workaround involves development of radiators based on electromagnetic (EM) simulations which enable accurate evaluation of performance, yet at high computational cost. A typical EM-driven design process requires a geometry template (either obtained based on engineering experience or derived from the literature). The stencil is then modified according to the needs and numerically tuned using the optimization algorithm of choice [14], [15]. The procedure can be realized iteratively in a setup where the experience-driven development of topology is followed by a sequence of (manual) geometry changes and performance-oriented optimizations in hope that the structure with a desirable performance will be obtained [16]-[18]. While proven useful, the outlined cognitive approach is subject to engineering bias manifested as a tendency to select topologies, and/or geometrical modifications known to the designer. Other bottlenecks of experience-driven design include substantial time expenditure and a limited number of requirements that can be (reliably) controlled in the course of topology development [19]. Furthermore, the process cannot guarantee that suitable solutions will be identified. From this perspective, unsupervised approaches that prioritize specification-oriented antenna generation over engineering insight appear to be promising alternatives.

Automatic antenna design involves direct adjustment of free-form topology according to the requirements. Adoption of the approach is subject to availability of generic EM simulation models, i.e., the ones that are not limited to specific pre-determined shapes, nor input parameters [20]-[27]. The generic topologies fall into two main categories that enable representation of radiator as: (i) a composition of primitives, or (ii) a set of coordinate points. For the first class, the geometries are encoded using binary matrices. Their contents correspond to activation/deactivation of partially overlapping (to ensure electrical connection) basic components (e.g., rectangles, or triangles) [20]-[24]. Having a one-to-one relation between the shape to be obtained and its matrix-based representation is convenient for maintaining topological consistency even for complex geometries with metallization inclusions. On the other hand, large matrices are required to ensure flexibility in terms of geometries that can be generated. Furthermore, the need to control both the distribution of primitives within the matrix and their physical (floating-point) dimensions reveals a mixed-integer nature of the synthesis problem which might necessitate adaptation of suitable optimization routines [28]. For the second class of methods, the point-based coordinates—interconnected using, e.g., splines, or line sections—are used to define the antenna outline. The scheme can be conveniently integrated with popular numerical algorithms [10], [25]-[27]. On the other hand, the appropriate distribution of points has to be maintained in the course of design optimization so as to prevent generation of infeasible (e.g., self-intersecting) solutions [9]. Moreover, similarly to the matrix-based models, ensuring suitable geometrical flexibility for coordinate-based shapes requires a large number of independent design parameters [10]. Another challenge, shared by both discussed approaches, is related to difficulties in determination of useful initial geometries for optimization.

The problem pertinent to identification of suitable design candidates (also referred to as starting points) for optimization can be mitigated, albeit to some extent, using population-based optimization (e.g., metaheuristics) [27], [29]. Capability of such algorithms to process a set of designs enables evaluation of diverse topologies at each iteration (and hence thorough exploration of the search space) [30]. On the other hand, due to the use of EM simulations to determine antenna performance, the outlined mechanism poses a serious challenge in terms of computational cost [10]. Consequently, practical applications of conventional metaheuristic algorithms are often limited to synthesis of simple structures represented using up to a few independent parameters and characterized by fairly standard responses [26], [27]. Unacceptable cost of population-based topology development can be reduced using surrogate-assisted methods [10], [31]. Their goal is to shift the burden associated with optimization from the expensive simulation model to a cheap data-driven surrogate identified using a set of training geometries and their EM-based responses [10], [19], [31]. While coupling surrogates with metaheuristic algorithms proved to be useful, their practical applicability often limited to low-dimensional problems [32]. Other methods involve extending the number of parameters that can be handled through identification of data-driven models within a previously confined region the search space [31], [61]. However, the approach is based on implicit assumption that the topology of the radiator and the parametric ranges required to identify useful fraction of the search space are known [31]. From this perspective, the outlined approaches are not suitable to automatic generation of multi-dimensional geometries.

Synthesis of antenna topologies in a local-search setup represents an interesting alternative to optimization using population-based methods. Owing to adjustment of only one design at a time, local algorithms require only a handful of EM simulations per iteration (to determine a descent direction). Furthermore, the design problem can be embedded into a surrogate-assisted framework in order to maintain manageable computational cost. A popular EM-driven design method involves adjustment of antenna geometries using a gradient-based routine embedded within a trust-region (TR) loop [10], [33], [34]. The method boils down to generation of a local data-driven model around the best available design solution. The model is then optimized and tuned or re-set based on the quality of the predicted solutions with respect to EM-based responses [34]. The outlined approach has been used for development of multi-parameter antennas, also in setups that implement variable-fidelity EM simulations or sequential adjustment of the design problem dimensionality [10]. Regardless of the demonstrated usefulness, successful local-search-based topology generation is subject to availability of suitable antenna geometry. In practice, the candidate designs are either obtained from empirical equations [35] or derived through approximation of shapes from the literature that are promising from perspective of the imposed design specifications [10]. Although local optimization of generic antenna models can be performed at an acceptable computational cost (especially when coupled with surrogate-assisted methods), the process is subject to limitations resulting from the mentioned cognition-driven determination of starting points. In summary, adaptation of the state-of-the-art techniques to automatic development of modern antennas characterized by non-standard performance responses is associated with issues pertinent to identification of appropriate topologies [10], [26], unacceptable design cost [27], [36], but also convergence of optimization algorithms [37]. From this perspective, the problem of unsupervised antenna design remains open.

In this work, a variable-fidelity framework for automatic, specification-oriented design of free-form antennas has been proposed. The procedure involves semi-random generation of topologies followed by their EM-based evaluation using a low-fidelity (i.e., coarsely discretized) model. For each design, the obtained responses are then embedded into a surrogate-assisted classification routine which enables evaluation of the antenna performance while re-scaling its dimensions at a negligible cost. The candidate solutions are accepted for local tuning if their responses do not violate the classifier-defined threshold. Numerical optimization of promising designs is performed using a gradient-based algorithm embedded within a TR framework. The resulting low-fidelity geometries are then used as starting points for TR-based fine-tuning (i.e., at the high-fidelity EM simulation model level), so as to obtain the final designs. The last step involves just a few iterations. The selected contributions of the work include: (i) development of a surrogate-assisted method for uniform scaling of antenna topologies, (ii) integration of the method into a robust design classification routine with a warm-start capability (i.e., re-use of already existing designs as potential candidates for new design problems), as well as (iii) implementation of the framework for automatic, cost-efficient generation of asymmetric antenna structures represented using multiple independent parameters. The proposed framework has been demonstrated based on six test cases concerning development of free-form patch antennas for IoT applications. The radiators are represented using over fifty design variables. Two sets of unconventional structures dedicated to operate within 5 GHz to 6 GHz and 6 GHz to 7 GHz frequency ranges are developed. The optimized antennas are characterized by operational bandwidths that vary from 17% to 20%, which vastly outperforms capabilities of conventional patch antennas that operate over a range of around 4%. The automatically generated (and optimized) structures are experimentally validated and compared against state-of-the-art radiators from the literature. Comparison of the proposed framework against the state-of-the-art algorithms, as well as evaluation of the obtained designs in terms of sensitivity to manufacturing tolerances are also considered.

## 2. Design Methodology

A variable-fidelity framework for automatic generation of free-form patch antennas is discussed here. To ensure that the work is self-contained, the formulation of the design problem is provided. Next, the mechanism for surrogate-assisted classification of quasi-randomly generated designs is explained. Then, the trust-region local optimization of the variable-fidelity EM-simulation model is discussed, and the section is concluded by the summary of the presented methodology.

*2.1 Problem Formulation*

Let $\boldsymbol{R}_f(\boldsymbol{x}) = \boldsymbol{R}_f(\boldsymbol{x}, f)$ be the high-fidelity response of the free-form antenna obtained over a frequency range of interest $f$ for the vector of input parameter $\boldsymbol{x}$. The design optimization task can be defined as [38]:

$$\boldsymbol{x}^* = \arg\min_{\boldsymbol{x} \in X} \left( U\left( \boldsymbol{R}_f(\boldsymbol{x}) \right) \right) \quad (1)$$

where $X$ is the feasible region of the search space, $\boldsymbol{x}^*$ is the optimum design to be found, and $U = U(\boldsymbol{R}_f(\boldsymbol{x}))$ is the objective function. EM-driven antenna optimization involves direct evaluation of $U$ starting from a known initial design $\boldsymbol{x}^{(0)} \in X$ based on the $\boldsymbol{R}_f$ model simulations [38].

In conventional design scenarios, both the topology of the radiator, as well as the starting point for (1) are determined based on the engineering experience to narrow down the search space to the region of interest and aid identification of $x^*$ using local routines [32], [37]. However, for unsupervised development of free-form radiators, the vector $x$ represents both, the structure shape and its dimensions. Hence, the problem of antenna generation involves determination of $X \subset X_n$ (where $X_n$ contains both the feasible and infeasible topologies) and identification of promising candidate geometries within $X$. For the sake of this study, the generation of antenna structures represented in the form of interconnected coordinates is assumed. Therefore, the set of requirements for the feasible topology includes lack of self-intersections within the generated shape, as well as electrical connection between the feed point and the driven element (i.e., a radiator defined using a set of points). Moreover, due to lack of data on the distribution of useful designs, the lower/upper bounds $l_b/u_b$ of $X$ should be relaxed [31]. Also, maintaining flexibility of geometries w.r.t. attainable performance is subject to their representation using multiple independent parameters. As a consequence, the search space is not only large (due to its high dimensionality) but also multi-modal which poses a substantial challenge in terms of its exploration using local routines.

Besides the challenges related to identification of useful starting point $x^{(0)}$, minimization of (1) is impractical due to a tremendous number of EM simulations required for the algorithm to converge when dealing with multi-dimensional problems. Instead, direct evaluations of $R_f$ can be replaced by an iterative procedure which gradually approximates the final design, through optimization of iteratively updated surrogate models $R_s^{(i)}$, $i = 0, 1, \ldots$, as [38]:

$$x^{(i+1)} = \arg\min_{x \in X}\left(U\left(R_s^{(i)}(x)\right)\right) \qquad (2)$$

where $x^{(i+1)}$ approximates the design $x^*$ and $R_s^{(i)}$ is obtained from EM simulations of the $R_f$ model performed at the candidate designs appropriately distributed around $x^{(0)}$. The computational cost of the process can be further reduced upon embedding (2) in a variable-fidelity framework [31], [39], where $R_s^{(i)}$ is first obtained from the underlying low-fidelity model $R_c$. The latter one is a simplified (hence inaccurate, yet computationally cheap) counterpart of $R_f$.

The problem considered in this work involves automatic generation of free-form antennas and their local tuning. The design is performed in two-steps that include: (i) generation and scaling of initial solutions, as well as (ii) their surrogate-assisted optimization in a variable-fidelity setup. The details concerning each step are given below.

*2.2 Surrogate-Assisted Generation of Initial Designs*

Determination of initial designs is performed using a quasi-random method. The coordinates that define the antenna outline are first generated, and then sorted so as to ensure lack of shape self-intersections (and hence maintain its feasibility). Next, the coordinates of feed are randomly generated while ensuring their confinement within the outline. Clearly, the use of random points distribution does not provide any control over the performance of the structure under development. Given that the typical design requirements involve determination of the acceptable reflection within the frequency of interest, the probability of obtaining a suitable starting point for multi-parameter design problem is low. On the other hand, complex geometry might offer a non-standard electrical response with multiple local minima spread over a relatively broad bandwidth. Figure 1 illustrates the response obtained for two randomly generated antennas along with the design requirements. Although the performance within the range of interest is poor for both structures, their responses feature resonances above the frequency of interest that could be leveraged for optimization. Given the relation between the antenna size and its electrical behavior, the geometry can be scaled to shift the desirable features of the response (e.g., resonances) to usable range. Hence, broadband evaluation and scaling of performance characteristics increases the change that randomly generated shapes will be useful for local optimization. The shift of antenna frequency response can be approximated, at a negligible cost, using regression [31], [38].

Let $R_{s.f}(x)$ be the surrogate obtained through shifting of the low-fidelity antenna model response using as $f_1 = \alpha f$, where $f$ is the original frequency sweep (cf. Section 2.1) and $\alpha$ is the scaling coefficient [38]. The construction and change of the $R_{s.f}$ responses is realized through approximation of $R_c(x, f)$ upon determination of a mapping between $\alpha$ and the antenna size [38]. The model is given as:

$$\alpha = \beta_0 c^2 + \beta_1 c + \beta_2 \qquad (3)$$

where $c$ represents the coefficient for uniform scaling of antenna dimensions and $\boldsymbol{\beta} = [\beta_0\ \beta_1\ \beta_2]^T$ are model weights. Let $X_q = \{x_{q.1}, \ldots, x_{q.t}\}$, $t = 1, \ldots, T$, be a set of quasi-randomly generated designs. Then let $y_{t.k} = c_k \cdot x_{q.t}$ and $R_{c.t.k} = R_c(y_{t.k})$ be the design scaled uniformly using $c_k = c_0 + \delta_k$ and its corresponding frequency response. The parameters $c_0$ and $\delta_k$ denote the nominal weighting coefficient and $k$th ($k = 0, 1, 2, \ldots, K$) element of the vector $\boldsymbol{\delta} = [0\ \delta_1\ \ldots\ \delta_k]^T$. The latter one contains the selected positive/negative components that alter the $c_0$ scale. The goal is to determine the relation between antenna size and frequency shift of its responses based on data extracted through evaluation of $T$ designs over $K$ weighting coefficients. For the obtained $R_{c.t.k}$ responses, the frequency of a specific resonance $r_{t.k}$ is selected and its shift with $c_k$ is then tracked. The result is a vector of design-specific frequency resonances $\boldsymbol{r}_t = [r_{t.0}\ r_{t.1}\ \ldots\ r_{t.K}]^T$ which is then normalized w.r.t. $r_{t.0}$ to obtain $\boldsymbol{\alpha}_t = r_{t.0}/\boldsymbol{r}_t$. The extraction is performed at all $T$ training designs. The resulting matrix of frequency scaling coefficients $A = [\boldsymbol{\alpha}_1\ \ldots\ \boldsymbol{\alpha}_t]^T$ is obtained for the given $C = \boldsymbol{1}c$ where $\boldsymbol{c} = [c_0\ \ldots\ c_K]$ and $\boldsymbol{1}$ is the column vector of ones. Finally, the weights $\boldsymbol{\beta}$ for (3) are obtained analytically using

the least squares method [40]. The algorithm can be summarized as follows:

1. Set $t = 1$, $k = 0$, $c_0$ and $\delta$;
2. Generate $y_{t,k}$ based on $x_{q,t}$, obtain $R_{c,t,k}$ based on EM simulations and extract $r_{t,k}$;
3. If $k = K$, calculate $a_t = r_{t,0}/r_t$ and go to Step 4; otherwise set $k = k + 1$ and go to Step 2;
4. If $t = T$, solve (3) using the least squares method, construct the $R_{s,f}$ model using the extracted $\beta$ and END; otherwise set $t = t + 1$ and go to Step 2.

It is worth noting that the components of $\beta$ can be also approximated analytically from $r_k = v/c_k$, where $v$ is the speed of light. However, the simulation-based coefficients provide a more realistic approximation of the frequency changes, especially when the selected $c_k$ weights are distant to $c_0$. Figure 2 shows a comparison of the scaling curves extracted based on EM simulations of multiple designs versus the analytical solution, as well as a response of the example antenna before and after scaling. Note that the considered $R_{s,f}$ model only adjusts the location resonances over frequency. Control of their level has not been considered to ensure a low cost of multi-dimensional designs generation (one $R_c$ simulation per solution). Although construction of the surrogate involves evaluation of a few $R_{c,t,k}$ responses, the process is performed only once for the given class of antennas. Furthermore, it can be integrated with the quasi-random generation of shapes so as to re-use the already available designs. Once the surrogate is constructed, each candidate solution is evaluated using the classification function. For a more detailed discussion on construction and evaluation of frequency-scaled surrogate models, see [31], [38].

*2.3    Surrogate-Assisted Classification of Designs*

Upon EM simulation, the responses obtained for quasi-randomly generated designs $x_q$ are tuned—based on the surrogate model $R_{s,f}(x_q) = R_{s,f}(x_q, c)$ responses—according to the given classification function to evaluate their usefulness as candidates for local optimization. The classifier is defined with respect to the acceptable threshold $E_t$ on the in-band response level. The optimization process involves identification of the scaling coefficient that provides the best in-band response of the generated model. It is defined as:

$$c^* = \arg\min_c \left( U_q \left( R_{s,f} \left( x_q, c \right) \right) \right) \tag{4}$$

where $c^*$ is the optimum scaling coefficient to be found and design objective $U_q$ is given as:

$$U_q = -E_t + \left\{ \max \left( R_{s,f} \left( x_q, c \right) \right) \right\}_{f_L \leq f \leq f_H} \tag{5}$$

Here, $f \in \boldsymbol{f}$ is the range of interest for the antenna response constrained by the lower and upper frequencies $f_L$ and $f_H$, respectively. It should be reiterated that minimization of (5) involves only one EM simulation of the $R_c(x_q)$ response that is fed to the surrogate of Section 2.2. The candidate design is accepted for optimization, if the final response obtained from (5) is below zero. The proposed cost-efficient classification mechanism can be summarized as follows:

1. Set $E_t$, $f_L$, and $f_H$;
2. Generate quasi-random candidate design $x_q$ and evaluate its $R_c(x_q)$ response;
3. Feed the $R_c(x_q)$ response to the surrogate model $R_s$ and solve (4) to identify $c^*$;
4. If $U_q(R_{s,f}(x_q, c^*)) \leq 0$, set $x^{(0)} = x_q$ and END; otherwise go to Step 2.

The design $x^{(0)}$ identified using the outlined procedure is then considered as the starting point for local optimization. It should be noted that, although the method does not guarantee identification of appropriate initial solution, implementation of uniform scaling into the generated topologies enables generation of promising designs for multi-dimensional problems at an acceptable computational cost (especially when compared to brute-force optimization using population-based metaheuristic methods) [26], [27], [30]. Furthermore, even if the obtained design candidate is of no use for the given problem, it is stored in a database for future re-use.

*2.4    Design Optimization Engine*

One of the main challenges related to optimization of multi-dimensional structures (regardless of the single- or multi-objective nature of the problem) is a large number of EM simulations required to evaluate the quality of solutions in the given iteration. Although the problem can be mitigated (to some extent) through parallelization of EM simulations during optimization, such an approach put a large demand on the computing hardware, as well as licenses for simulation software [31], [41]. Instead, the problem can be embedded in a trust-region (TR) optimization loop where a set of $j = 0, 1, \ldots$ approximations to the final solution $x^*$ is generated using the iteratively updated surrogate model. The TR-based design task is formulated as follows [34]:

$$x^{(j+1)} = \arg\min_{\|x-x^{(j)}\|\le\lambda} \left(U\left(G^{(j)}(x)\right)\right) \qquad (6)$$

The parameter $\lambda$ represents a TR radius (here, initial value of $\lambda = 1$ is used), and $G$ is a local linear model constructed from the EM simulation model $R$ responses. Note that, depending on the antenna design stage (cf. Section 2.5), the surrogate $G$ is generated either based on $R = R_c$, or $R = R_f$ data, respectively. The model is given as:

$$G^{(j)}(x) = R(x^{(j)}) + J(x^{(j)})(x - x^{(j)}) \qquad (7)$$

The Jacobian $J$ is computed using a large-step central finite difference (FD), and is given as [42]:

$$J(x^{(p)}) = \left[\ldots \; \frac{1}{h_d^{(p)}}\left(R\left(x^{(p)} + \frac{h_d^{(p)}}{2}\right) - R\left(x^{(p)} - \frac{h_d^{(p)}}{2}\right)\right) \; \ldots\right]^T \qquad (8)$$

where $h_d^{(p)} = \sigma \cdot x_d^{(p)}$ represents the perturbation with respect to $d$th ($d = 1, \ldots, D$) dimension of $x^{(p)} = [x_1^{(p)} \ldots x_d^{(p)} \ldots x_D^{(p)}]^T$ design, $\sigma$ is the step size scaling coefficient, and $h_d^{(p)} = [0 \ldots h_d^{(p)} \ldots 0]^T$, respectively. Note that minimization of (6) can be performed in a mode where the Jacobian (8) is re-set after each successful iteration which corresponds to setting $p = j$. Alternatively, a static Jacobian defined around the initial design (i.e., with $p = 0$) can be used to reduce the cost of (8) while narrowing-down the final design from a good starting point, e.g., when a precise tuning of high-fidelity designs is considered. Moreover, the perturbation coefficient $\sigma$ can be adjusted to alter exploitation capabilities of the algorithm.

The routine (6) is controlled based on coefficient $\rho = (U(R(x^{(j+1)})) - U(R(x^{(j)})))/(U(G^{(j)}(x^{(j+1)})) - U(G^{(j)}(x^{(j)})))$, which represents a ratio between the obtained and predicted improvement of the objective function $U$ response. When $\rho > 0.75$, the TR radius is set to $\lambda = 2\lambda$. Alternatively, $\lambda = \lambda/3$ is set when $\rho < 0.25$. The algorithm is terminated when $\lambda^{(j+1)} < \varepsilon$ (here, $\varepsilon = 10^{-2}$), $\|x^{(j+1)} - x^{(j)}\| < \varepsilon$, or $\|U(R(x^{(j+1)})) - U(R(x^{(j)}))\| < \varepsilon$. It should be noted that the computational cost of (6) corresponds to $2D + 1$ EM simulations for each successful iteration in a mode where re-construction of $G$ is considered, i.e., when $p = j$ in (8). Additional simulation is required for each unsuccessful design, or step that does not involve $G$ model re-set ($p = 0$). For more comprehensive discussion on TR-based optimization, see [33], [34], [37].

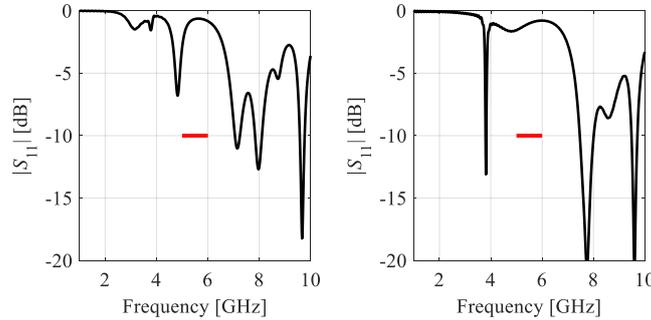

Fig. 1. Reflection responses obtained for two randomly generated multi-dimensional antenna topologies along with the example design specifications (red line). Note that the responses above the frequency range of interest could be leveraged for local optimization upon suitable scaling.

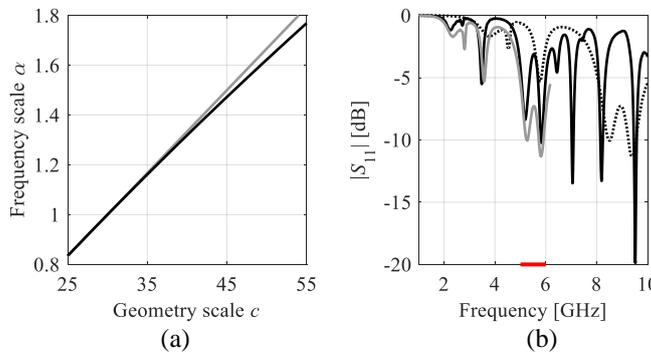

Fig. 2. Surrogate-assisted geometry tuning: (a) scaling curves based on analytical form (gray) and extracted from EM simulations (black), as well as (b) comparison of the $R_c$ at $c_0$ (…) against $R_{s.f}$ (gray) and $R_c$ response (black) for a random design scaled w.r.t. the frequency range of interest (red line). Note that the surrogate is only capable of shifting the resonances to ensure low cost of scaling.

*2.5    Design Framework Summary*

The proposed algorithm integrates: (i) quasi-random generation of designs, their cost-efficient uniform scaling and classification, as well as (ii) surrogate-assisted optimization of promising candidate solutions embedded in a variable-fidelity setup. The latter stage retains control not only over the quality of EM simulations used for generation of the local surrogate, but also the need for its re-set and the step size scaling (here, $\sigma = 0.02$ has been selected based on a series of numerical experiments) for the FD perturbations used in (8). Assuming availability of the generic antenna model, the presented framework can be summarized as follows (see Fig. 3 for a block diagram):

1. Set $D$, $f$, $f_L$, $f_H$, $E_t$ and identify the $R_{s,f}$ model using the method of Section 2.2;
2. Generate/evaluate the candidate design $x_q$ using algorithm of Section 2.3;
3. Determine $x^{(0)}$ using routine of Section 2.3;
4. Set $R = R_c$, $\sigma$ and obtain $x_{tmp}^*$ by solving (6) with Jacobian re-set ($p = j$);
5. If $R(x_{tmp}^*)$ violates the design requirements, set $\sigma = \sigma/2$, $x^{(0)} = x_{tmp}^*$ and re-optimize design using (6) to obtain $x_c^*$; otherwise set $x_c^* = x_{tmp}^*$;
6. Set $j = 0$, $R = R_f$, $\sigma$, $x^{(j)} = x_c^*$ and obtain $x_f^*$ by solving (6) with static Jacobian ($p = 0$), END.

Owing to reliability of the Jacobians determined using central FD, the need to re-set the optimization in Step 4 is uncommon. The mechanism involves optimization of the $R_c$ model in a TR loop in order to identify the high-quality $x_c^*$ design followed by re-optimization of the structure based on $R_f$ model responses. The goal of the process is to narrow-down the promising region of the search space using coarse EM simulations, and its further exploitation using a limited number of accurate (yet expensive) fine model evaluations. It should be emphasized that the use of static Jacobian in Step 6 stems from a high quality of the identified $x_c^*$ solution. Consequently, the number of EM simulations required to obtain $x_f^*$ is low and the algorithm is terminated after two to three iterations (hence, identification of $J(x^{(0)})$ has the highest contribution to the overall cost of the last design step). For a more comprehensive discussion on variable-fidelity optimization of antenna structures see [38], [39].

### 3.    Numerical Results

In this section, the proposed framework for automatic generation and optimization of free-form antennas is demonstrated through optimization of a planar structure represented using a swarm of coordinate points [27]. Two case studies spanning across a total of six geometries are considered. The design problems involve unsupervised development of unconventional planar radiators for IoT applications dedicated to operate within the frequency bands of: (i) 5 GHz to 6 GHz and (ii) 6 GHz to 7 GHz, respectively. For both test problems the number of independent design parameters is set to $D = 53$ (including the coefficient required for uniform scaling of topologies). The corner frequencies for EM model simulations are set to a range from 1 GHz to 10 GHz. The objective function for optimization of the identified candidate designs involves minimization the antenna reflection within the bands of interest and it is given as:

$$U = \frac{1}{M}\sum_{m=1}^{M} \max\left(R_m(x) - R_{\max}, 0\right)^2 \qquad (9)$$

Here, $R_m(x) \in R(x)$, $m = 1, .., M$ represents the antenna reflection at the $m$th frequency point within the $f_L \leq f_m \leq f_H$ band of interest (cf. Section 2.3), whereas $R_{\max}$ is the target reflection level for the optimization process (here, $R_{\max} = -11$ dB is selected to ensure a slight margin w.r.t. the target $R_{goal} = -10$ dB). The reference scaling factor for antenna geometry and the threshold for the classifier are set to $c_0 = 30$ mm and $E_t = -5$ dB, respectively.

*3.1    Generic Antenna Model*

Unsupervised specification-oriented design of free-form antennas is subject to availability of a generic EM-simulation model. Consider a universal planar radiator shown in Fig. 4 [29]. The structure is implemented on a Rogers AD255C substrate ($h = 1.524$ mm, $\varepsilon_r = 2.55$, $\tan\delta = 0.0013$). It consists of a shape- and dimension-independent patch represented using a swarm of points defined in a cylindrical coordinate system. Geometric feasibility of the radiator (i.e., lack of outline of self-intersections) is ensured by the definition of each angular coordinate with respect to the previous one and normalization of their cumulative sum to a range from 0 to $2\pi$. It should be emphasized that (contrary to e.g., [27], [35]) the symmetry planes of the model are intentionally not defined to enhance the flexibility of the radiator in terms of attainable responses while enabling generation of non-standard radiation patterns (e.g., with dual-lobe or off-centered maxima). The patch is fed through a cylindrical probe. It should be reiterated that coordinates of the probe are generated randomly unless their confinement within the generated outline is achieved.

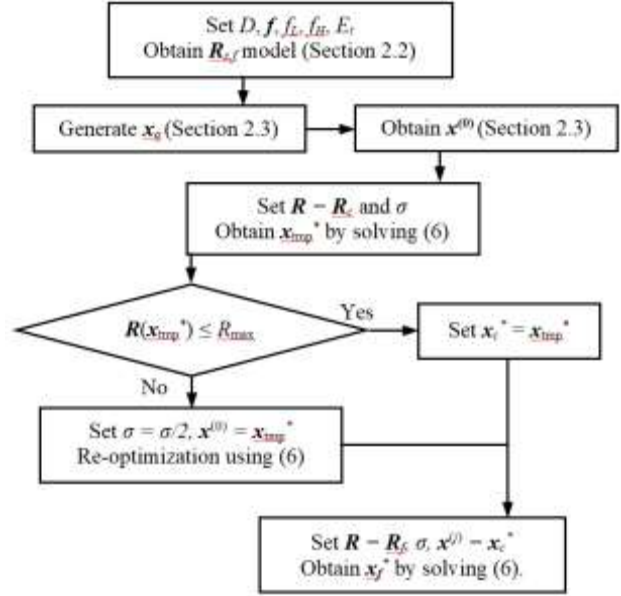

Fig. 3. A block diagram of the proposed optimization framework.

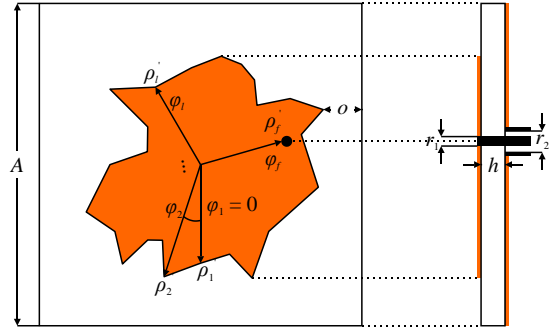

Fig. 4. Generic antenna EM simulation model used for optimization with highlight on design parameters.

The antenna design parameters are defined using a vector $x = [c\ \rho_f\ \varphi_f\ \rho\ \varphi]^T$, where $\rho = [\rho_1\ \rho_2\ \ldots\ \rho_l]^T$ and $\varphi = [\varphi_1\ \varphi_2\ \ldots\ \varphi_l]^T$ are the radial and angular coordinates used to represent the patch outline ($l = 1, \ldots, L$), whereas $\rho_f$ and $\varphi_f$ define the position of the probe feed, and $c$ is the antenna scaling factor (cf. Section 2.2). The fixed parameters include offset from the patch to the substrate edge $o = 5$ mm, as well as concentric feed radii of $r_1 = 0.615$ mm and $r_2 = 1.415$ mm (to ensure a 50 Ohm input impedance), respectively. Note that $c$ is defined in mm, whereas the remaining components of $x$ are unit-less. The overall number of adjustable parameters used to represent the structure is $D = 2L + 3$ (hence, $L = 0.5(D - 3)$). The substrate is a square with an edge given as $A = A_1 + 2o$, whereas the patch size is $A_1 = 2c \cdot \max(\rho)$. The topology is considered feasible within the following lower and upper limits $l_b = [c_{0.l}\ 0\ \varphi_{0.l}\ 0.1I\ 0.01I]^T$ and $u_b = [c_{0.h}\ \rho_{0.f}\ \varphi_{0.h}\ 0.9I\ 0.8I]^T$, where $I$ is an $L$-dimensional vector of ones. Note that the parameters $c_{0.l} = c^{(0)} - 2$, $c_{0.h} = c^{(0)} + 3$, $\rho_{0.f} = \max(\rho^{(0)})$, $\varphi_{0.l} = \varphi^{(0)} - \pi / 2$, and $\varphi_{0.h} = \varphi^{(0)} + 3\pi / 2$ are determined based on the starting point $x^{(0)}$ obtained as explained in Section 2.3.

The antenna model is implemented in CST Microwave Studio and evaluated using its time domain solver [43]. In order to support the variable-fidelity optimization, two models of the structure have been implemented. The low-fidelity model $R_c$ is represented using a coarsely discretized mesh comprising ~121,000 hexahedral cells. Other simplifications include the use of a lossless substrate, representation of the metallization as an infinitely thin perfect conductor, and relaxed criteria for time-domain algorithm convergence [44]. Its average evaluation time on the AMD EPYC 7282 system with 32 GB RAM is 60 s. The high-fidelity model $R_f$ is discretized using ~ 410,000 cells (simulation cost: 110 s) and introduces loss to dielectric substrate, as well as non-zero metallization [44].

*3.2 First Case Study*

The first design problem concerns specification-oriented development of antenna topologies dedicated to operate within $f_L = 5$ GHz to $f_H = 6$ GHz frequency range of interest. The process has been performed using the framework outlined in Section 2.5. The first step involves development of a mapping between the multiplicative frequency coefficient and the uniform scaling factor (3) for the surrogate model $R_{s.f}$ of Section 2.2. To obtain the vector of coefficients $\beta$, a set of $T = 3$ quasi-random designs $Y_q = \{y_{q.1}, y_{q.2}, y_{q.3}\}$ has been generated along with the following vector of re-scaling coefficients $\boldsymbol{\delta} = [0\ –5\ 15]^T$. The latter corresponds to $\boldsymbol{c} = [30\ 25\ 45]^T$, respectively. Then, a combination of the designs from $\boldsymbol{X}_q$ and scales from $\boldsymbol{c}$ has been evaluated using the simulation model of Section 3.1. The resulting responses have been processed to extract the $\boldsymbol{r}_t$ resonances and the vector $\boldsymbol{\beta} = [–10^{–4}\ 0.037\ –0.047]^T$ has been extracted using the least squares method. The model has been used along with classification mechanism of Section 2.3 to determine the promising candidate designs for local optimization. A total of three test designs, i.e., $\boldsymbol{x}_1^{(0)}$, $\boldsymbol{x}_2^{(0)}$, and $\boldsymbol{x}_3^{(0)}$ have been generated (see Appendix *A* for details). Note that for each starting point, an average change of the coefficient *c* by almost 12 mm (w.r.t. the reference value of $c_0 = 30$) has been introduced. Furthermore, neither of the non-scaled quasi-random design was considered acceptable by the classifier for use as the candidate solution which demonstrates usefulness of the proposed frequency-adjusting $\boldsymbol{R}_{s.f}$ surrogate.

The obtained design candidates have been then optimized using the TR engine of Section 2.4. For $\boldsymbol{x}_1^{(0)}$, the low-fidelity topology $\boldsymbol{x}_{c.1}^*$ has been found after 15 iterations of the algorithm (6). In the next step, the design has been used as a starting point for the TR-based tuning that embeds high-fidelity EM simulations (with a static Jacobian). The final design $\boldsymbol{x}_{f.1}^*$ has been obtained after just 2 iterations (for dimensions, see Appendix *A*). Figure 5 illustrates the shape and reflection responses of the antenna topology at each stage of the variable-fidelity optimization process. The obtained results indicate that, upon topology generation, the optimization based on $\boldsymbol{R}_c$ evaluations is suitable for identification of the search space region that contains the optimum design. The high-fidelity model with a static Jacobian is used only to perform a fine-tuning of the response. Hence, its computational cost is low (given the limited number of required EM simulations). Despite the relatively minor modifications of the topology shape in Fig. 5(c), the change of the maximum in-band response corresponds to a few dB, which indicates the significance of proper—i.e., high-fidelity-based—antenna tuning. It is worth noting that $\boldsymbol{R}_c$ response of the initial design in Fig. 5(b) violates the specified $E_t$ threshold for the classification function. As already explained in Section 2.2, this is because the $\boldsymbol{R}_{s.f}$ surrogate only supports frequency adjustment and thus no control over the response level is maintained. At the same time, the scaled $\boldsymbol{R}_c$ response is distorted due to unconventional shape of the patch. Notwithstanding, a comparison of the $\boldsymbol{R}_c(\boldsymbol{x}_1^{(0)})$ and $\boldsymbol{R}_{s.f}(\boldsymbol{x}_1^{(0)})$ demonstrates that shifting the resonances (see Fig. 2(b)), increases usefulness of the design compared to other, non-scaled solutions. Another important remark is that lower reflection at $c_0$ compared to *c* suggests that, upon re-scaling, minimization of the response within the range of interest should be achievable. This might not be the case when the response features closely located in frequency, yet topologically unrelated resonances (cf. Fig. 1(a)).

The overall numerical cost for the $\boldsymbol{x}_{f.1}^*$ design generation corresponds to 666.6 $\boldsymbol{R}_f$ model simulations (~20.4 h of CPU-time using a single machine). The cost comprises 6, 93, and 916 $\boldsymbol{R}_c$ evaluations for identification of the model (3), quasi-random identification of the initial design $\boldsymbol{x}_{c.1}^{(0)}$, and its TR-based optimization, as well as 113 $\boldsymbol{R}_f$ simulations for fine-tuning of the final design. Note that extraction of $\boldsymbol{\beta}$ for (3) is performed only once as coefficients can be re-used for scaling w.r.t. different frequency ranges of interest. The computational cost can be considered low, especially given a large number of design parameters used to represent the antenna topology.

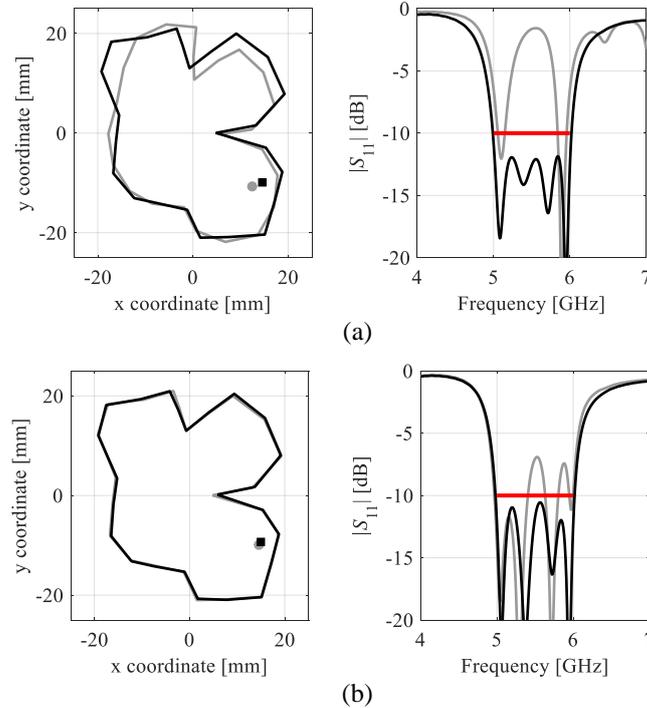

Fig. 5. Unconventional antenna developed using the proposed framework – the geometry and frequency responses of the structure represented using: (a) $\boldsymbol{R}_c$ model at $\boldsymbol{x}_1^{(0)}$ (gray), and $\boldsymbol{x}_{c.1}^*$ (black) designs, as well as (b) $\boldsymbol{R}_f$ model at $\boldsymbol{x}_{c.1}^*$ (gray) and $\boldsymbol{x}_{f.1}^*$ (black) designs (for dimensions, see Appendix A).

Table 1: First Case Study – Comparison Of The Optimized Designs

| Design | $A_1$ [mm] | $f_L$ [GHz] | $f_H$ [GHz] | BW [GHz] | BW [%] | Gain [dBi] @ $f$ [GHz] 5 | 5.5 | 6 |
|---|---|---|---|---|---|---|---|---|
| $x_{f.1}^*$ | 42.5 | 4.99 | 6.02 | 1.03 | 18.7 | 7.49 | 2.70 | 5.18 |
| $x_{f.2}^*$ | 44.7 | 4.96 | 6.07 | 1.10 | 19.8 | 2.61 | 8.10 | 5.65 |
| $x_{f.3}^*$ | 44.6 | 4.99 | 6.02 | 1.03 | 18.7 | 3.05 | 7.18 | 6.30 |

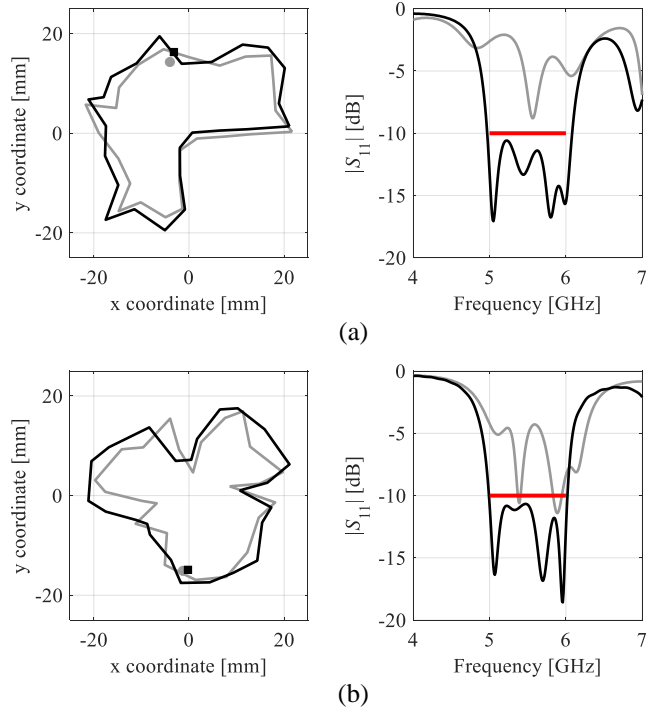

Fig. 6. Topologies dedicated for 5GHz to 6 GHz bandwidth – geometries at the initial (gray) and optimized (black) designs, as well as their corresponding reflection responses obtained for antenna: (a) 2 and (b) 3.

For the remaining two designs, the topologies and responses at $x_2^{(0)}$, and $x_{f.1}^*$, as well as $x_3^{(0)}$ and $x_{f.3}^*$ are compared in Fig. 6 (their dimensions are given in Appendix A). It is worth noting that $x_{f.1}^*$ and $x_{f.3}^*$ geometries exhibit a noticeable similarity of the outlines, but also substantially different locations of feeds. The computational cost of designs generation amount to 784.3 $R_f$ simulations (~23.9 h of CPU-time) for $x_{f.2}^*$ and 910.2 $R_f$ evaluations (~27.8 h) for $x_{f.3}^*$, respectively. The costs are low given high dimensionality of antenna model.

Table 1 compares the generated antennas in terms of bandwidth (BW), maximum gain, and dimensions. The structures exhibit broadband responses of around 20%, and their maximum gain varies from ~3 dBi to ~8 dBi. It should be emphasized that the structures have only been optimized according to the electrical performance requirements. Consequently, their specific gains are obtained merely as a by-product of the reflection-oriented design.

The generated radiators $x_{f.1}^*$, $x_{f.2}^*$, and $x_{f.3}^*$ can be enclosed within squares with the side length $A_1$ (cf. Section 3.1) of 42.2 mm (footprint: ~1781 mm$^2$), 44.7 mm, and 44.7 mm (size: ~1998 mm$^2$), respectively. At the same time, a rectangular radiator implemented on the same substrate and dedicated to operate at 5.5 GHz frequency features a side length of only 15.7 mm and overall size of ~246 mm$^2$ [13]. However, the smaller antenna offers only 3.6% bandwidth at the $R_{goal} = -10$ dB level (cf. Section 3). Increasing the size of the radiator to 31.4 mm × 33.3 mm (footprint: ~1046 mm$^2$)—so as to shift its second harmonic to around 5.5 GHz—results in BW enhancement to 6.7%. At the same time, the design $x_{f.1}^*$ features operational range of 18.7%, which is 5.2 and 2.8 times broader compared to conventional patches, yet at the expense of 7.2 and 1.7 times larger footprint. A comparison of the discussed antenna topologies and their responses at the optimized designs is shown in Fig. 7. It should be noted that, although the response of the larger rectangular patch features local minima around 6.2 GHz, it cannot be used to further enhance the structure BW due to insufficient number of degrees of freedom.

*3.3 Second Case Study*

The second case study concerns the development of antennas dedicated to operate within the $f_L = 6$ GHz to $f_H = 7$ GHz bandwidth. The $R_{s.f}$ surrogate used for re-scaling of quasi-random designs remains the same as in Section 3.2. It should be reiterated that all of the candidate solutions obtained for the first case study have been stored in the database. The available designs have been re-used for classification (cf. Section 2.3). Two candidates, i.e., $x_4^{(0)}$ and $x_5^{(0)}$ have been based on warm-start solutions. Hence, the selection of promising design solutions has been performed at a cost of only two $R_c$ simulations (for evaluation of the accepted, scaled geometries). Moreover, to further demonstrate the capability of repurposing the initial points for various performance requirements, the design $x_6^{(0)}$ has been obtained by re-scaling $x_2^{(0)}$ according to the specified operational bandwidth (with the resulting size-adjusting coefficient of $c$ = 37.64 mm). Next, the geometries have been optimized using the framework of Section 2.5.

The geometries and reflection responses of the antennas at $x_4^{(0)}$, $x_5^{(0)}$, $x_6^{(0)}$, $x_{f.4}^*$, $x_{f.5}^*$, and $x_{f.6}^*$ designs are shown in Fig. 8 (for dimensions, see Appendix B), whereas a qualitative comparison of their performance is provided in Table 2. Similarly as in Section 3.2, relatively small (yet noticeable) topology changes between the initial and optimized designs have substantial effects on antennas electrical performance. Each of the optimized topologies features a relatively broadband behavior (with bandwidth at $R_{goal}$ exceeding 1 GHz, or around 17% w.r.t. the center frequency of 6.5 GHz). At the same time the side length of squares that confine the patches span from 33 mm (for $x_{f.5}^*$) to 35.75 mm (for $x_{f.4}^*$), respectively.

As indicated in Table 2, the maximum gain of radiators (again, obtained as a by-product of reflection-oriented optimization) varies from 6.7 dBi for $x_{f.4}^*$ to 8 dBi for $x_{f.6}^*$. It is worth noting that the maximum gains for $x_{f.2}^*$ and $x_{f.6}^*$ differ only by 0.1 dB, which implies a significant contribution of antenna shape to the overall radiation performance and hence the potential to obtain a non-standard radiation behavior.

The computational cost of generating the optimized designs amount to 586.8 (~17.9 h of CPU-time), 581.8 (~17.8 h), and 637.7 $R_f$ simulations (~19.5 h) for $x_{f.4}^*$, $x_{f.5}^*$, and $x_{f.6}^*$, respectively. Again, the cost is low having in mind high dimensionality of the geometries (53 variables per design). It should be reiterated that all of the initial design candidates have been obtained by re-purposing the quasi-random solutions generated for the problem of Section 3.2 (i.e., the ones already stored in the database). The acceptable performance of the optimized solutions (all of them fulfill the design requirements) show that available quasi-random geometries can be used to warm-start the proposed topology-generation procedure. Given that a large database of design candidates is available, the cost of selecting appropriate structure can be greatly reduced compared to identification of new quasi-random solution for each test problem. As demonstrated, the cost has been reduced from an average of 388 $R_c$ simulations (for the first case study) to merely one per design.

## 4. Measurements

The optimized antennas have been fabricated and measured. Photographs of the selected prototypes are shown in Fig. 9, whereas comparisons of EM simulations and measurements in terms of reflection (for $x_{f.1}^*$ to $x_{f.6}^*$ designs) are provided in Fig. 10. The discrepancy between maxima of the in-band EM-based and experiment-based responses vary from 0.02 dB for $x_{f.5}^*$ to 1.42 dB for $x_{f.2}^*$. Furthermore, the measured reflection characteristics for $x_{f.2}^*$ and $x_{f.4}^*$ designs slightly violate the requirement concerning $R_{goal}$ = –10 dB. Overall, the agreement between the simulated and experimentally obtained electrical performance characteristics is acceptable.

The far-field responses of manufactured antennas have been evaluated in a non-anechoic regime using in-house developed measurement system and corrected using appropriate post-processing algorithm [45], [62]. Comparisons of the simulated and measured radiation patterns obtained in yz-plane (cf. Fig. 9(b)) at 5.5 GHz for $x_{f.1}^*$, $x_{f.2}^*$, $x_{f.3}^*$, and at 6.5 GHz for $x_{f.4}^*$, $x_{f.5}^*$, and $x_{f.6}^*$ designs, respectively, are shown in Fig. 11. It should be reiterated that the obtained characteristics are by-products of the reflection-oriented optimization. Most of the structures exhibit the non-standard dual-lobe responses with local maxima at around 35° and –40°, respectively. In this regard, the $x_{f.3}^*$ design is an outlier with only one pronounced radiation maxima at 5°. Again, the resemblance between the obtained responses is acceptable, especially given that measurements have been performed in non-standard propagation conditions. It is worth noting that slight yet noticeable discrepancies between the simulations and measurements can be attributed to manufacturing tolerances (as small changes of antenna shape might notably affect its performance, see Fig. 5), simplified EM models (i.e., rudimentary representation of concentric feeds, lack of solder joints, etc.), manual assembly of structures, as well as challenging measurement conditions [45]. For more detailed discussion on the utilized test equipment and post-processing routine, see [45], [62].

Table 2: Second Case Study – Comparison of the Optimized Designs

| Design | $A_1$ [mm] | $f_L$ [GHz] | $f_H$ [GHz] | BW | | Gain [dBi] @ $f$ [GHz] | | |
|---|---|---|---|---|---|---|---|---|
| | | | | [GHz] | [%] | 6 | 6.5 | 7 |
| $x_{f.4}^*$ | 35.75 | 5.99 | 7.18 | 1.19 | 18 | 2.52 | 6.70 | 2.67 |
| $x_{f.5}^*$ | 32.99 | 5.98 | 7.04 | 1.06 | 16 | 6.48 | 6.73 | 2.59 |
| $x_{f.6}^*$ | 34.80 | 5.96 | 7.09 | 1.13 | 17 | 4.16 | 8.00 | 6.96 |

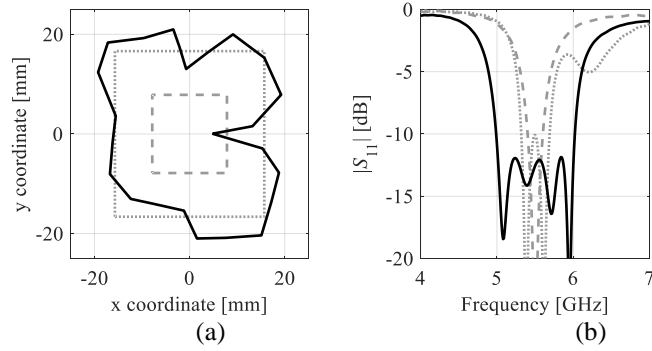

Fig. 7. Comparison of standard (– –) and enlarged (···) rectangular patches against the free-form antenna $x_{f.1}^*$ (—) in terms of: (a) size and shape (in-scale), as well as (b) in-band responses. Note that the enlarged patch shifts harmonic frequency to around 5.5 GHz and features additional resonance.

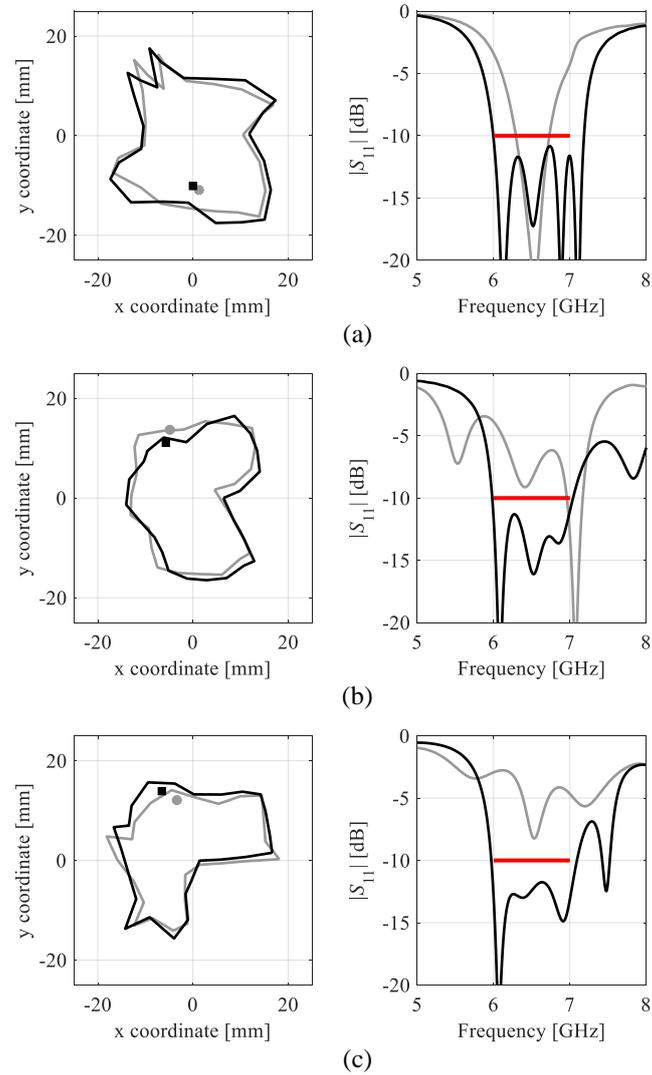

Fig. 8. Topologies dedicated for 6 GHz to 7 GHz bandwidth – geometries at initial (gray) and optimized (black) designs and their corresponding reflection responses obtained for antenna: (a) 4, (b) 5, and (c) 6.

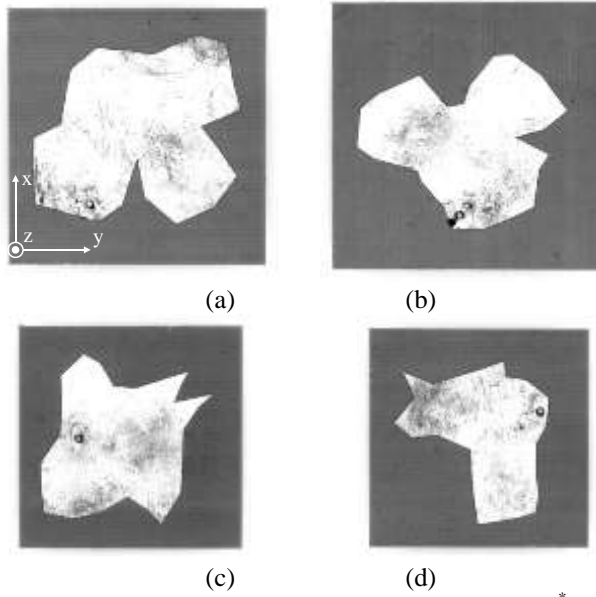

(a) (b) (c) (d)

Fig. 9. Photographs (in-scale) of the manufactured antenna prototypes: (a) $x_{f.1}^*$, (b) $x_{f.3}^*$, (c) $x_{f.4}^*$, and (d) $x_{f.6}^*$. Overlay axes denote orientation of the antennas for the radiation pattern measurements.

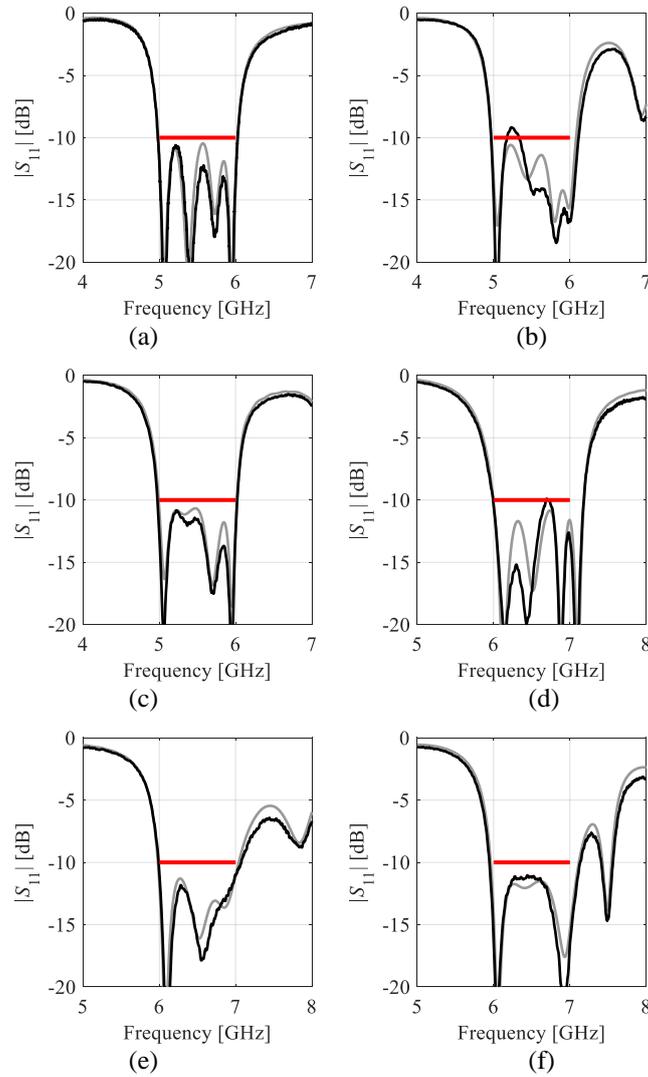

Fig. 10. In-band reflection – comparison of EM simulations (gray) and measurements (black) for designs: (a) $x_{f.1}^*$, (b) $x_{f.2}^*$, (c) $x_{f.3}^*$, (d) $x_{f.4}^*$, (e) $x_{f.5}^*$, and (f) $x_{f.6}^*$.

## 5. Comparisons and Discussion

In this section, the proposed method has been benchmarked against competitive algorithms. Furthermore, sensitivity of the antenna prototypes to fabrication tolerances and their comparison against radiators from the literature have been performed. Potential applications of the optimized designs in IoT systems, challenges related to application of the framework to automatic development of antennas, as well as benefits resulting from warm-start of topology generation have also been discussed.

*5.1 Benchmark of the Method*

The proposed framework has been compared against other algorithms in terms of computational cost and performance. The considered methods include TR-based routine of Section 2.4, as well as a standard population-based evolutionary algorithm with mating restrictions, tournament-selection, and elitism (population size: $2D$, maximum number of iterations: 250, probability of mutation and crossover: 0.5 and 0.2) [46]. Specific benchmarks involve: (i) TR-based optimization without surrogate-assisted adjustment of antenna size (cf. Section 2.3) [9], [10], [34], as well as metaheuristic-based optimization in setups (ii) with scaling of designs from the initial population according to the classification function and (iii) without scaling (with $c$ = 30 mm) [20], [26], [27]. For the sake of fair comparison, the algorithm (i) has been executed starting from a modified design $x_1^{(0)}$ (with $c$ = 30 mm), whereas initial population for the evolutionary routine comprises the 106 quasi-random designs (twice the number of design parameters, i.e., $2D$) already stored in the database and includes the $x_1^{(0)}$ solution (cf. Appendix A). Furthermore, feasibility has been ensured for all of the designs generated during optimization using appropriate pre-conditioning mechanism. For all of the test cases, minimization of the objective function (9) has been considered. The lower and upper bounds on the frequency of interest are the same as in Section 3.2.

The results collected in Table 3 indicate that, for the considered test scenarios, the presented framework not only generates the best solution, but also at the lowest computational cost. Inferior $U(x^*)$ responses obtained using methods (i) and (iii) suggest that the scaling mechanism of Section 2.3 is useful for improving the quality of initial solutions. It is worth noting, that for (iii), the algorithm has been terminated after 100 iterations due to lack of objective function improvement. Interestingly, the method (i) generated better design than the population-based routines which indicates difficulties related to reliable exploration of complex and highly-dimensional search spaces using limited populations, as well as the significance of having good starting point. Furthermore, the cost of metaheuristic optimization is an order of magnitude higher compared to the local-based methods and corresponds to over a week of computations (using a single-CPU machine). Consequently, the considered population-based methods seem to be of limited use for the problem at hand.

*5.2 Statistical Analysis*

Regardless of high resemblance between EM simulations and measurements (demonstrated in Section 4), the design $x_{f.2}^*$ slightly violates the in-band $R_{goal}$ = –10 dB threshold (cf. Section 3). The requirement has not been fulfilled despite determination of target level at $R_{max}$ = –11 dB to ensure a slight margin for degradation of antenna performance. As already mentioned, the discrepancy might be affected by the fabrication tolerances. The latter ones can be quantified as a result of statistical analysis (also referred to as yield estimation) based on a Monte Carlo (MC) simulations [47], [48]. The problem involves evaluation of the responses for a set of $M$ random perturbations generated—according to the assumed probability distribution—around the nominal (i.e., optimized) design. The manufacturing yield represents a ratio of the designs that fulfill the specifications to all of the available perturbations around the given $x_f^*$ solution. It can be estimated as [47]:

$$Y = \frac{1}{M} \sum_{m=1}^{M} H\left(U_1\left(x^{(m)}\right)\right) \tag{10}$$

where $x^{(m)}$ is $m$th perturbed design. The function $H$ is equal to 1 when the $U_1(x) \leq 0$ and 0 otherwise. The objective function considered here is of the form:

$$U_1(x) = R_{goal} - \max\{R(x)\}_{f_L \leq f \leq f_H} \tag{11}$$

It should be noted that direct EM-driven MC simulations are prohibitively expensive as they require tremendous number of perturbed designs to ensure satisfactory accuracy. Here, surrogate-assisted yield estimation has been performed for all designs [47], [48]. The considered tests include evaluation of 10000 perturbations—generated with uniform (maximum deviation $\sigma_{max}$ = 0.05 mm) and Gaussian (mean $\mu$ = 0 and standard deviation $\sigma$ = 0.03 mm) probability distributions— around the optimized designs. The surrogates for each test case have been generated at a cost corresponding to $2D$ high-fidelity simulations (a total of 624 $R_f$ evaluations).

The results gathered in Table 4 demonstrate that the manufacturing yield calculated w.r.t. (11) is close to 100% for all of the considered antennas except for $x_{f.4}^*$. Notwithstanding, the target level $R_{max}$ has not been achieved for the latter one. Furthermore, high yield for $x_{f.2}^*$ suggest that violation of the performance requirement is not due to the manufacturing tolerances, but other effects. It is expected that manual assembly (installation, soldering, etc.) of the probe feed at the edge of the antenna outline (cf. Fig. 8(b)) is the likely cause of the discrepancy. Figure 12 shows MC simulations performed w.r.t. the $x_{f.2}^*$ and $x_{f.4}^*$ designs using a few hundred randomly generated samples (Gaussian distribution). It should be noted that deviations of responses from the nominal solutions are low. For more comprehensive discussion on surrogate-assisted MC analysis, see [47], [48].

Table 3: Benchmark of the Algorithm

| Method | No of EM simulations | | Total cost | | $U(x^*)$ |
| --- | --- | --- | --- | --- | --- |
| | $R_c$ | $R_f$ | $R_f$ | [h] | |
| (i) | 1271 | – | 693.3 | 21.2 | 13.9 |
| (ii) | 26500 | – | 14455 | 442 | 16.9 |
| (iii) | 10600 | – | 5782 | 177 | 25.2 |
| This work | 1015 | 113 | 666.6 | 20.4 | 0.01 |

Table 4: Optimized Antennas – Fabrication Yield

| Distribution | | Nominal designs | | | | | |
| --- | --- | --- | --- | --- | --- | --- | --- |
| | | $x_{f.1}^*$ | $x_{f.2}^*$ | $x_{f.3}^*$ | $x_{f.4}^*$ | $x_{f.5}^*$ | $x_{f.6}^*$ |
| Gaussian | $Y$ [%] | 99.0 | 99.8 | 96.7 | 75.5 | 99.9 | 100 |
| Uniform | $Y$ [%] | 99.3 | 99.7 | 95.9 | 76.2 | 100 | 100 |

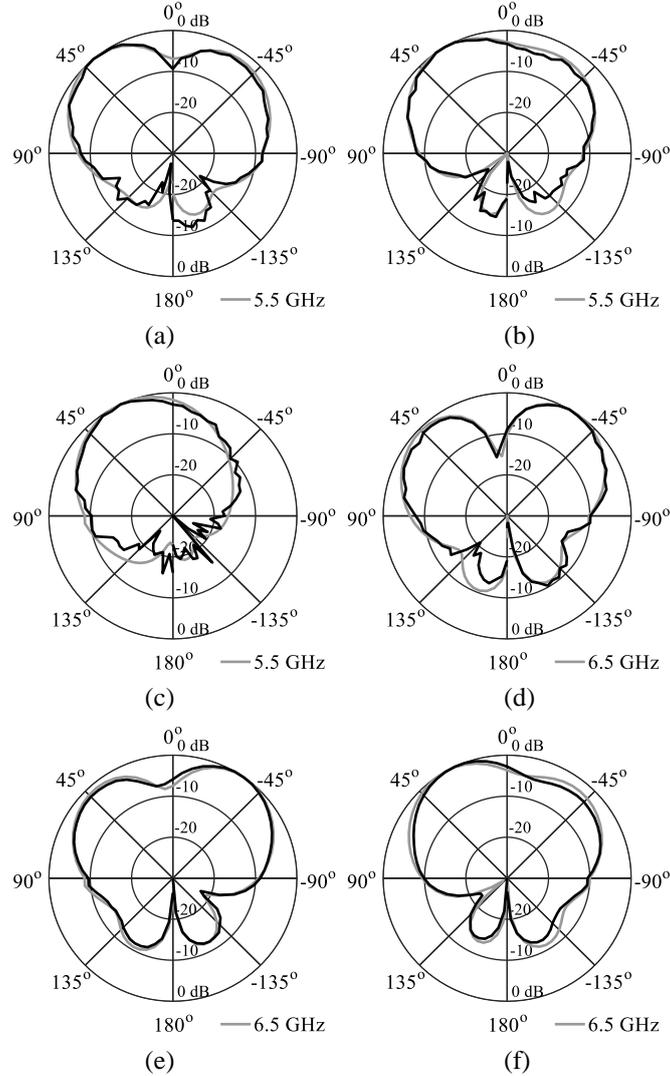

Fig. 11. Radiation patterns in yz-plane (cf. Fig. 9(b)) – comparison of EM simulations (gray) and measurements (black) for the designs: (a) $x_{f.1}^*$, (b) $x_{f.2}^*$, (c) $x_{f.3}^*$, (d) $x_{f.4}^*$, (e) $x_{f.5}^*$, and (f) $x_{f.6}^*$, respectively.

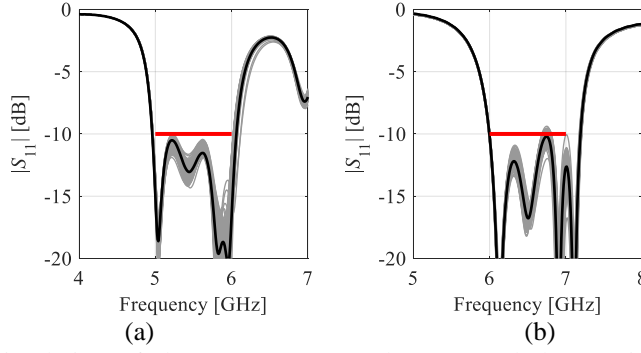

Fig. 12. Monte Carlo simulation of the antenna (gray plots) around the nominal design (black): (a) $x_{f.2}^*$ (estimated yield: 100%) and (b) $x_{f.4}^*$ (yield: 68.9%). Red bars denote the desired operational bandwidth.

Table 5: Benchmark Of The Optimized Designs

| Geometry | Footprint | | $f_0$ | BW | | Gain | $A \times B$ |
|---|---|---|---|---|---|---|---|
| | $[\lambda_g \times \lambda_g]$ | $[\lambda_g^2]$ | [GHz] | [GHz] | [%] | [dBi] | [mm × mm] |
| Rectangular[#] | 0.84 × 0.89 | 0.74 | 5.50 | 0.37 | 6.72 | 4.54 | 31.4 × 33.3 |
| [51] | 0.40 × 0.60 | 0.24 | 5.90 | 0.27 | 4.58 | 6.02 | 11.1 × 16.7 |
| [9] | 0.46 × 0.65 | 0.30 | 5.50 | 0.87 | 15.8 | – | 13.6 × 19.5 |
| [50] | 0.55 × 0.57 | 0.31 | 2.46 | 0.18 | 7.30 | 2.90 | 37.5 × 38.7 |
| [52] | 0.49 × 0.65 | 0.32 | 6.70 | 0.35 | 5.23 | 1.54 | 12.0 × 16.0 |
| [53] | 0.59 × 0.85 | 0.50 | 3.86 | 0.08 | 2.17 | 5.34 | 25.2 × 36.0 |
| [54] | 0.93 × 0.93 | 0.87 | 5.30 | 0.31 | 5.83 | 1.48 | 28.9 × 28.9 |
| [55] | 0.79 × 2.08 | 1.64 | 7.52 | 2.11 | 28.0 | 8.00 | 19.0 × 50.0 |
| [56] | 1.36 × 1.36 | 1.85 | 5.40 | 0.17 | 3.05 | 7.90 | 54.0 × 54.0 |
| $x_{f.3}^*$ | 1.19 × 1.19 | 1.41 | 5.50 | 1.03 | 18.7 | 7.49 | 44.6 × 44.6 |
| $x_{f.5}^*$ | 1.04 × 1.04 | 1.08 | 6.50 | 1.06 | 16.3 | 6.73 | 33.0 × 33.0 |

[#] Obtained for the enlarged patch of Section 3.2

### 5.3 Comparison Against Other Antennas and Discussion

Two antennas developed using the proposed framework—i.e., $x_{f.3}^*$ and $x_{f.5}^*$ designs—have been compared in terms of performance and size against the state-of-the-art radiators from the literature [9], [51]-[56]. The considered figures of interest include BW, maximum gain and size. For fair comparison, the dimensions of each antenna have been expressed in terms of guided wavelength $\lambda_g$ calculated w.r.t. the electrical parameters of the substrate used structure for implementation, as well as its center frequency $f_0$. Moreover, the footprints of unconventional radiators have been expressed as the area of the smallest $A \times B$ rectangle that encloses the patch. The results collected in Table 5 indicate that all but one benchmark antenna offer BW below 1 GHz. Although the structure [55] features bandwidth of over 2 GHz (28%) and high gain, it has been implemented on two-fold thicker substrate (compared to $x_{f.3}^*$ and $x_{f.5}^*$) and incorporates multi-layer topology with a few radiating components. The antenna of [9] offers BW of 0.88 GHz (i.e., 16%), yet its gain has not been reported. Overall, most of the considered benchmark structures are characterized by inferior performance (both in terms of bandwidth and gain), but also smaller electrical dimensions. Having that in mind, the presented automatically developed radiators maintain acceptable compromise between the physical dimensions and performance characteristics.

Owing to exceptionally broad bandwidth (as for patch-based antennas) and dual-lobe behavior, the radiators generated using the proposed framework might find application for optimization of IoT systems in terms of balancing their cost and performance. A good example would include their integration as a part of heterogenic nodes for in-door positioning systems [57], [58]. Mentioned nodes are normally optimized to work with omnidirectional antennas which are useful for operation in large open areas, yet non-optimal for environments such as office buildings, or other sites with complex arrangements. Transition to heterogenic architecture tailored to propagation environment might reduce the overall cost of system installation and maintenance while ensuring its acceptable performance [59]. Another potential application might include cost-efficient adaptation of wireless connectivity gear (e.g., using WiFi 6) to the operational conditions with emphasis on maintaining low number of access points required to ensure coverage [60].

### 5.4 Challenges

The challenges related to application of the proposed design framework stem from the multi-dimensional nature of the problem. They include: (i) unpredictable cost associated with quasi-random generation of (scaled) topologies and (ii) the risk of getting stuck in poor local optimum.

Stochastic nature of the proposed topology generation process makes *a priori* assessment of its cost impossible. Furthermore, the method does not guarantee that acceptable geometry will be found. On the other hand, for the considered test cases, promising (upon scaling) starting designs have been generated at an average cost being below 400 $R_c$ model simulations (cf. Section 3.3). However, warm-start of the design identification procedure based on the quasi-random solutions already available in the database affects the estimation. At the same time, none of the candidate designs was considered acceptable when scaling was neglected. Clearly, the results and analyses do not imply that algorithmic topology-development cannot be performed using the geometries neglected by the classifier of Section 2.3. Nonetheless, when the optimization is performed using the local-search algorithm (as in this work), a good starting point increases the chance of convergence to the acceptable design solution (cf. Section 5.1). It is worth noting that, given availability of broadband responses, the designs gathered in the database can be deemed universal, as scaling (or other screening methods) can be applied to identify promising candidates for specifications that include broadband (demonstrated), multiband, or other behavior of interest. It is worth noting that dimensionality of the designs (and their responses) from the database can also be interpolated in order to adjust the number of independent design parameters to requirements, or capabilities of the optimization algorithm [10], [35].

The TR engine considered in this work is a robust optimization method. However, its local-search nature poses the risk of the algorithm to get stuck in a poor optima. As already indicated, the latter is partially associated to the high dimensionality and multi-modal character of the considered automatic topology development problem. On the other hand, the results of Section 5.1 demonstrate that identification of useful solutions is also a challenge for global optimization routines [26], [27]. Furthermore, the mentioned algorithms require tremendous number of EM simulations to converge (here, over a week of computations), which challenges their feasibility for direct development of multi-dimensional structures. At the same time, various implementations TR-based methods proved to be useful for optimization of antenna structures characterized by many independent variables [33], [37]. Other methods such as sequential optimization of components followed by gradual increase of the number of parameters have also been considered. However, they seem to negatively affect the granularity of the search space [10], [35]. When compared against metaheuristics, TR methods seem to represent an acceptable trade-off between the computational cost and the quality of achievable solutions.

The discussed challenges indicate potential directions for future research on reliable and automatic specification-oriented design of free-form antennas. Regardless, it has been demonstrated that the proposed framework is capable of generating intricate structures featuring performance specifications that are beyond the reach for conventional patch-based radiators.

## 6. Conclusion

In this work, a framework for specification-oriented surrogate-assisted development of free-form antennas has been proposed. The method involves automatic generation and cost-efficient adjustment of quasi-random topologies, as well as their evaluation using the classification function. The promising design candidates are then optimized within a TR framework. The process is performed in a variable-fidelity setup, where the final design is first approximated based on low-fidelity EM simulations and then fine-tuned using the accurate EM simulation model. Free-form geometries are represented using a generic EM simulation model. The proposed approach has been demonstrated based on two case studies concerning design of broadband radiators (a total of six antenna structures) dedicated to operate within the 5 GHz to 6 GHz and 6 GHz to 7 GHz bands. The bandwidths of the obtained topologies span from 16% to 20%, respectively. The average computational cost of structure development corresponds to 695 high-fidelity simulations (~21.2 h of CPU-time) per design. Experimental validation of manufactured prototypes indicates a good resemblance between the simulations and measurements, which confirms the validity of the proposed universal EM simulation model. Most of the generated structures feature dual-lobe radiation characteristics which make them of potential use for IoT applications that include development of heterogenic in-door localization systems, wireless connectivity networks dedicated to complex propagation environments, or point-to-multi-point communication. The proposed framework has been benchmarked against competitive algorithms. Furthermore, the optimized antennas have been evaluated in terms of sensitivity to manufacturing tolerances and favorably compared against the state-of-the-art radiators from the literature. The obtained results demonstrate usefulness of the proposed framework for specification-oriented development of antennas at an acceptable computational cost.

Future work will focus on incorporation of nature-inspired design generation methods to increase diversity of geometries, as well as incorporation of inclusions within the antenna outlines while reducing the cost associated with determination of useful topologies. Other directions of research will include application of the proposed methodology for generation of volumetric structures, incorporation of yield maximization into the design framework, as well as development of free-form radiators while accounting for multiple performance figures (such as in-band reflection, axial ratio, or gain) at a time.

## Acknowledgement
This work was supported in part by the National Science Centre of Poland Grant 2021/43/B/ST7/01856.

## Appendix
The appendix contains initial and optimized designs of the free-form antenna structures obtained using the proposed design framework.

## A. First Case Study – Design Parameters

The vectors of adjustable parameters that represent the free-form antennas generated and accepted for optimization using the routines of Section 2 are listed below:

1. $x_1^{(0)}$ = [42 0.47 5.72 0.24 0.4 0.53 0.54 0.52 0.42 0.28 0.38 0.52 0.52 0.49 0.4 0.32 0.32 0.33 0.35 0.35 0.36 0.48 0.58 0.64 0.62 0.56 0.45 0.24 0 0.05 0.21 0.32 0.3 0.11 0.18 0.12 0.05 0.29 0.32 0.28 0.41 0.49 0.46 0.45 0.37 0.46 0.09 0.23 0.25 0.25 0.24 0.2 0.16]$^T$;
2. $x_2^{(0)}$ = [45 0.31 1.58 0.57 0.5 0.53 0.58 0.48 0.37 0.34 0.37 0.34 0.31 0.28 0.42 0.34 0.27 0.3 0.43 0.34 0.38 0.35 0.22 0.09 0.13 0.25 0.37 0.57 0 0.19 0.25 0.19 0.14 0.15 0.41 0.34 0.38 0.41 0.29 0.1 0.31 0.44 0.43 0.09 0.18 0.36 0.22 0.03 0.28 0.87 0.15 0.06 0.03]$^T$;
3. $x_3^{(0)}$ = [45 0.38 4.63 0.19 0.33 0.44 0.33 0.42 0.35 0.2 0.06 0.17 0.32 0.28 0.36 0.44 0.39 0.23 0.17 0.3 0.23 0.37 0.42 0.44 0.4 0.35 0.41 0.19 0 0.04 0.11 0.42 0.37 0.14 0.23 0.01 0.55 0 0.6 0.2 0.42 0.12 0.22 0.18 0.13 0.52 0.15 0.39 0.33 0.34 0.42 0.23 0.18]$^T$.

The vectors obtained from the $x_1^{(0)}$ design after TR-based optimization of the low- and high-fidelity (starting from $x_{c.1}^*$) EM models:

1. $x_{c.1}^*$ = [41.94 0.5 5.78 0.21 0.41 0.58 0.59 0.57 0.43 0.32 0.42 0.5 0.48 0.54 0.47 0.29 0.3 0.36 0.37 0.35 0.37 0.52 0.57 0.66 0.6 0.57 0.46 0.25 0.12 0.09 0.25 0.35 0.35 0.14 0.21 0.14 0.08 0.32 0.35 0.28 0.39 0.45 0.44 0.45 0.33 0.44 0.08 0.26 0.26 0.26 0.25 0.19 0.16]$^T$
2. $x_{f.1}^*$ = [41.96 0.5 5.81 0.22 0.41 0.57 0.59 0.57 0.43 0.31 0.41 0.5 0.48 0.54 0.47 0.29 0.3 0.37 0.38 0.36 0.37 0.51 0.57 0.66 0.59 0.56 0.46 0.21 0.01 0.09 0.24 0.33 0.33 0.14 0.21 0.15 0.1 0.3 0.33 0.26 0.37 0.42 0.41 0.42 0.31 0.42 0.09 0.24 0.24 0.25 0.24 0.18 0.16]$^T$

The design parameters that represent the free-form antennas optimized from $x_2^{(0)}$ and $x_3^{(0)}$ using the framework of Section 2.5:

1. $x_{f.2}^*$ = [45.03 0.33 1.48 0.57 0.53 0.61 0.59 0.51 0.36 0.29 0.4 0.32 0.34 0.32 0.39 0.29 0.32 0.35 0.51 0.4 0.47 0.38 0.23 0.11 0.12 0.25 0.38 0.57 0 0.19 0.25 0.2 0.16 0.14 0.38 0.33 0.38 0.4 0.29 0.1 0.31 0.43 0.43 0.11 0.22 0.37 0.25 0.05 0.29 0.77 0.16 0.04 0.04]$^T$;
2. $x_{f.3}^*$ = [44.84 0.35 4.67 0.23 0.36 0.48 0.45 0.44 0.4 0.24 0.15 0.15 0.34 0.43 0.48 0.48 0.41 0.29 0.25 0.27 0.33 0.42 0.43 0.42 0.45 0.37 0.39 0.24 0.01 0.1 0.16 0.41 0.36 0.18 0.23 0.11 0.5 0.13 0.49 0.2 0.37 0.15 0.24 0.17 0.17 0.47 0.18 0.32 0.28 0.29 0.38 0.21 0.2]$^T$.

## B. Second Case Study – Design Parameters

The design vectors accepted for optimization using the algorithm of Section 2.3:

1. $x_4^{(0)}$ = [42 0.27 4.88 0.27 0.37 0.44 0.34 0.27 0.26 0.36 0.41 0.26 0.36 0.45 0.32 0.23 0.23 0.37 0.43 0.35 0.36 0.35 0.39 0.44 0.52 0.46 0.37 0.27 0 0.26 0.07 0.36 0.42 0.53 0.3 0.03 0.15 0.1 0.02 0.18 0.48 0.48 0.08 0.14 0.36 0.36 0.39 0.34 0.26 0.15 0.23 0.26 0.35]$^T$;
2. $x_5^{(0)}$ = [35.99 0.36 1.94 0.13 0.25 0.37 0.45 0.49 0.43 0.39 0.34 0.37 0.43 0.41 0.33 0.35 0.38 0.33 0.4 0.48 0.47 0.47 0.51 0.48 0.49 0.37 0.25 0.13 0 0.13 0.19 0.25 0.20 0.27 0.33 0.34 0.33 0.32 0.15 0.35 0.36 0.3 0.34 0.24 0.19 0.27 0.27 0.25 0.28 0.13 0.12 0.16 0.52]$^T$;
3. $x_6^{(0)}$ = [37.64 0.31 1.58 0.57 0.5 0.53 0.58 0.48 0.37 0.34 0.37 0.34 0.31 0.28 0.42 0.34 0.27 0.3 0.43 0.34 0.38 0.35 0.22 0.09 0.13 0.25 0.37 0.57 0 0.19 0.25 0.19 0.14 0.15 0.41 0.34 0.38 0.41 0.29 0.1 0.31 0.44 0.43 0.09 0.18 0.36 0.22 0.03 0.28 0.87 0.15 0.06 0.03]$^T$.

The optimized vectors obtained using the framework of Section 2.5:

1. $x_{f.4}^*$ = [42.37 0.25 4.77 0.29 0.38 0.45 0.37 0.28 0.27 0.37 0.45 0.27 0.35 0.42 0.34 0.23 0.25 0.38 0.45 0.44 0.35 0.33 0.44 0.49 0.55 0.48 0.39 0.29 0 0.27 0.08 0.38 0.42 0.52 0.28 0.06 0.17 0.14 0.05 0.16 0.43 0.44 0.08 0.13 0.34 0.34 0.37 0.32 0.26 0.19 0.25 0.25 0.33]$^T$;
2. $x_{f.5}^*$ = [35.89 0.36 2.09 0.16 0.26 0.4 0.45 0.48 0.51 0.42 0.32 0.39 0.39 0.37 0.4 0.41 0.35 0.33 0.36 0.44 0.45 0.46 0.48 0.47 0.49 0.37 0.26 0.16 0 0.15 0.22 0.26 0.21 0.26 0.3 0.33 0.33 0.31 0.16 0.33 0.35 0.3 0.34 0.25 0.2 0.26 0.25 0.25 0.26 0.13 0.11 0.17 0.51]$^T$;
3. $x_{f.6}^*$ = [37.71 0.34 1.82 0.53 0.53 0.54 0.56 0.49 0.4 0.32 0.37 0.41 0.36 0.3 0.38 0.31 0.29 0.33 0.5 0.38 0.46 0.36 0.22 0.17 0.13 0.25 0.39 0.52 0.02 0.18 0.25 0.17 0.14 0.17 0.39 0.33 0.38 0.39 0.29 0.13 0.31 0.45 0.45 0.13 0.21 0.39 0.21 0.05 0.28 0.77 0.18 0.09 0.06]$^T$.